\documentclass[12pt]{amsart}
\usepackage{amssymb}
\usepackage{geometry}

\theoremstyle{definition}

\theoremstyle{remark}

\numberwithin{equation}{section}

\input{tcilatex}
\begin{document}
\title{On the K-Ring of the Classifying Space of the Dihedral Group}
\author{Mehmet K\i rdar}
\address{Department of Mathematics, Faculty of Arts and Science, Nam\i k
Kemal University, Tekirda\u{g}, Turkey}
\email{mkirdar@nku.edu.tr}
\date{9 May, 2012}

\begin{abstract}
We describe the K-ring of the classifying space of the dihedral group in
terms of generators and the minimal set of relations by emphasising the
connection with the polynomials arising in the KO-rings of lens spaces and
demonstrating the generators of the filtrations of the Atiyah-Hirzebruch
Spectral Squence.
\end{abstract}

\maketitle







\section{Introduction}

Let $D_{2n}$ denote dihedral group with $2n$ elements. In this note, by
collecting enough evidence from the the Atiyah-Hirzebruch Spectral Sequence
(AHSS) and essentially by making use of the structure of the representation
ring $R(D_{2n}),$ we will describe $K(BD_{2n})$ in terms of generators and
the minimal set of relations by means of the Atiyah-Segal Completion Theorem
(ASCT).

The $K$-rings of the $BD_{2n}$ are studied classically in [3]. For the
computations in the extended cohomology theory and its connective version,
see [1]. We note that the zeroth part of the\ $K$-cohomology and hence the
relations here should be extracted from the results in [1]. But, they are
not explicitly written there.

First we deal with the case where $n$ is odd and we will have a complete
description.

Then, we will consider the even case which turns out to be quite complicated
compared to the odd case. Thanks to [4], we could be able to find the
minimal relations in that case quickly. We also explain the twisted nature
of the AHSS by demonstrating the generators of its filtrations.

\section{The Odd Case}

Let $n=2k+1$ and $k\geq 1.$ Isomorphism classes of the irreducible complex
representations of $D_{2n}$ are the two one dimensional representations $%
\eta ,1$ and the $k$ two dimensional representations $\rho _{i}$ where $%
1\leq i\leq k.$ Actually, assuming that $\rho _{0}=1+\eta ,$ for any $i\in Z$
the symbol $\rho _{i}$ makes sense and we have the relations $\eta
^{2}=1,\rho _{i}=\eta \rho _{i}=\rho _{n-i}=\rho _{n+i}$ and $\rho _{i}\rho
_{j}=\rho _{i+j}+\rho _{j-i}$ for all $i,j\in Z$ in the representation ring $%
R(D_{2n})$.

Let us denote the complex vector bundles induced from these representations
over $BD_{2n}$ by the same letters.

The integral group cohomology of $D_{2n},$ for $n$ odd, is given in [2] as%
\begin{equation*}
H^{p}(BD_{2n};Z)=\left\{ 
\begin{tabular}{l}
$Z$ \\ 
$Z_{2}$ \\ 
$Z_{n}\oplus Z_{2}$ \\ 
$0$%
\end{tabular}%
\right. 
\begin{tabular}{l}
$p=0$ \\ 
if $p=4s+2$ \\ 
if $p=4s,$ $s\geq 1$ \\ 
if $p$ is odd.%
\end{tabular}%
\end{equation*}%
and since the cohomology in odd dimensions is zero, AHSS collapses on the
second page and we can read the filtrations of the ring $K(BG)$ from the
cohomology.

We define the reductions $v=\eta -1$ and $\phi =\rho _{1}-2$ in the ring $%
K(BD_{2n}).$ Since the representations $\eta $ and $\rho _{1}$ generate the
representation ring $R(D_{2n})$, $v$ and $\phi $ generate $K(BD_{2n})$ due
to ASCT. The relation $\eta ^{2}=1$ in our new variables is 
\begin{equation*}
v^{2}=-2v.
\end{equation*}%
From $\rho _{1}=\eta \rho _{1}$ we get $(\phi +2)v=0$ or equivalently 
\begin{equation*}
\phi v=v^{2}.
\end{equation*}

The classifying space $BZ_{2}$ is the infinite dimensional real projective
space and its $K$-ring is generated by $v$ modulo the relation $v^{2}=-2v$
where $v=\eta -1$ is the reduction of the tautological complex line bundle $%
\eta $ over $BZ_{2}.$

We have two natural group injections of $Z_{2}$ in the group $D_{2n}.$ Let
us pick one of them $Z_{2}<D_{2n}$. This induces a ring homomorphism $%
i:K(BD_{n})\rightarrow K(BZ_{2})$. Now we should have $i(v)=v$ by naturality
of the constructions. Since $v^{s},$ $s\geq 1,$ is generating $E_{\infty
}^{2s,-2s}$ filtration of $K(BZ_{2})$ in AHSS, the same is true for $%
K(BD_{n})$ by naturality and this explains the $Z_{2}$ parts of the
cohomology and the filtrations $E_{\infty }^{2s,-2s}=E_{2}^{2s,-2s}$ of the
AHSS. Now we should explain the remaining $Z_{n}$ parts on dimensions $4s$.

The classifying space $BZ_{n}$, for $n\geq 3$, is an infinite dimensional
lens space and its $K$-ring is generated by $\mu $ modulo the relation $%
(1+\mu )^{n}=1$ where $\mu =\eta -1$ is the reduction of the tautological
complex line bundle $\eta $ over $BZ_{n}.$ We have a natural group injection 
$Z_{n}<D_{2n}$ which induces a homomorphism $j:K(BD_{n})\rightarrow
K(BZ_{n}) $, as explained above for the naturally identified $Z_{2}$
subgroups of the dihedral group. Then we have $j(v)=0$ and $j(\phi )=\eta
+\eta ^{-1}-2$ in $K(BZ_{n}).$

Now, this last reduced bundle is the generator $c(w)=\eta +\eta ^{-1}-2$ of
the image of the $KO$-ring of $Z_{n}$ under the complexification $%
c:KO(BZ_{n})\rightarrow K(BZ_{n})$ which is injective. Here $w$ is the
realification of $\mu $ and it generates $KO(BZ_{n}).$ Recall that when $n$
is odd 
\begin{equation*}
KO(BZ_{n})=Z\left[ w\right] /(wf_{n}(w))
\end{equation*}%
where 
\begin{equation*}
wf_{n}(w)=\psi ^{\frac{n+1}{2}}(w)-\psi ^{\frac{n-1}{2}}(w)
\end{equation*}
and where $\psi ^{i}(\cdot )$ are the polynomials defined by $\psi
^{i}(w)=\dsum\limits_{j=1}^{i}\frac{\binom{i}{j}\binom{i+j-1}{j}}{\binom{2j-1%
}{j}}w^{j}$ for $i\geq 1$. Note that this is how the Adams operations are
acting on $w$ and $KO(BZ_{n})$ is just the fixed set of the complex
conjugation.

The polynomial $f_{n}(w)$ is explicitly given by%
\begin{equation*}
f_{n}(w)=n+\sum_{j=1}^{\frac{n-3}{2}}\frac{%
n(n^{2}-1^{2})(n^{2}-3^{2})...(n^{2}-(2j-1)^{2})}{2^{2j}.(2j+1)!}w^{j}+w^{%
\frac{n-1}{2}}.
\end{equation*}%
We will call $f_{n}(w)$ the minimal polynomial of the ring $KO(BZ_{n}).$ If $%
n=p$ is an odd prime number then $f_{n}(\cdot )$ is the minimal polynomial
of the number $2\cos (\frac{2\pi }{p})-2$. This can be proved by using the
Eisenstein's criterion and factorization of odd indexed Chebishev
polynomials.

Now, since $E_{\infty }^{4s,-4s}$ filtrations of $K(BZ_{n})$ and $K(BD_{n})$
coincides on $Z_{n}$ parts by means of $j$, and since $\phi v=v^{2}=-2v,$ we
deduce that $\phi $ satisfy a relation in the form $\phi f_{n}(\phi )+cv=0$
where $c$ is an integer to be determined. Multiplying this relation by $v$
and by using $\phi v=v^{2}$ we have $v^{2}(f_{n}(v)+c)=0$ and by using $%
v^{2}=-2v,$ we conclude that $c=-f_{n}(-2).$

We did a little calculation by using the various relations, essentially the
iterative relation $\rho _{i}\rho _{j}=\rho _{i+j}+\rho _{j-i}$ of the
representation ring, with some help of the compute option of the Scientific
WorkPlace. We observed that, indeed, for small values of $k$, $\phi $
satisfies a relation in the form $\phi f_{n}(\phi )\pm v=0$. Let's explain
why $c$ is $\pm 1.$ It is known that $f_{n}(\phi )$ is the square root of
the polynomial 
\begin{equation*}
\frac{2T_{n}(\frac{\phi +2}{2})-2}{\phi }
\end{equation*}%
where $T_{n}(\cdot )$ is the Chebishev polynomial of degree $n$ and $n$ is
odd. This last expression is just equal to $1$ when we substitute $\phi =-2,$
because $T_{n}(0)=0$. We conclude that $f_{n}(-2)=\pm 1.$

All these information fit together very nicely except the filtration
absurdity in the main relation%
\begin{equation*}
\phi f_{n}(\phi )=f_{n}(-2)v.
\end{equation*}%
On the left hand side we have something which seems to be of filtration at
least four $4$ while the right hand side is of filtration $2$. This may be
modified by claiming that $Z_{n}$ part of the filtration $E_{\infty
}^{4s,-4s}$ is generated by $\phi -v.$ In fact, this argument seems to be
justified in [3]. This kind of a generator in a filtration can be seen as
some kind of twisting phenomenon in the $K$-cohomology. But in any case, we
have the following

\bigskip

\textbf{Theorem 1: }Let $n\geq 3$ be odd. Then 
\begin{equation*}
K(BD_{2n})=Z[v,\phi ]\diagup (v^{2}+2v,v\phi +2v,\phi f_{n}(\phi
)-f_{n}(-2)v)
\end{equation*}%
where $f_{n}(\cdot )$ is the minimal polynomial of $KO(BZ_{n})$; in
particular the minimal polynomial of $2\cos \frac{2\pi }{n}-2$ when $n$ is
prime.$\blacksquare $

Let us make some concluding remarks. It is interesting that the even part of
the integral cohomology of the group, and thus the whole integral
cohomology, is completely detected by and detects the $K$-ring when $n$ is
odd. Since, the representations of the dihedral groups are real, the rings $%
K(BD_{2n})$ and $KO(BD_{2n})$ should be isomorphic as rings, but,
geometrically they carry different information. Their spectral sequences are
different and $KO(BD_{2n})$ should be related to some of the cohomology with 
$Z_{2}$ coefficients.

And note also that the curious constant in the theorem above seems to be%
\begin{equation*}
f_{n}(-2)=\sin \frac{k\pi }{2}+\cos \frac{k\pi }{2}
\end{equation*}%
where $n=2k+1.$

\section{The Even Case}

The even case is quite complicated with respect to the odd case. But,
recently, Kirdar and Sevilay, [4], described the $K$-ring of the generalized
quaternion groups by a minimal set of relations and we realized that the
relations of the representation ring here are very similar to the relations
there. And amazingly, we could be able to easily deduce the minimal set of
relations for the dihedral groups, the even case, with some little
modifications.

Since $D_{4}$ is $Z_{2}\times Z_{2},$ the space $BD_{4}$ is the product of
two infinite dimensional real projective spaces, and the $K$-ring can be
described by Atiyah's Kunneth Theorem for $K$-cohomology.

So, let us take $n=2k$ and $k\geq 2$ from now on. The representation ring $%
R(D_{4k})$ is generated by the four one dimensional representations $1,\eta
_{1},\eta _{2}$ and $\eta _{3}=\eta _{1}\eta _{2}$ and the $k-1$ two
dimensional representations $\rho _{i}$ where $1\leq i\leq k-1.$ Actually, $%
\rho _{i}$ make sense for any integer $i$, due to the standard definition of
these representations by matrices.

The various relations in $R(D_{4k})$ are: The start $\rho _{0}=1+\eta _{3}$
and the end $\rho _{k}=\eta _{1}+\eta _{2};$ $\eta _{1}^{2}=\eta
_{2}^{2}=\eta _{3}^{2}=1;$ $\eta _{r}\rho _{i}=\rho _{k-i}=\rho _{k+i}$ for $%
r=1,2$ and $\eta \rho _{i}=\rho _{i};$ $\rho _{n+i}=\rho _{i};$ and the main
relation $\rho _{i}\rho _{j}=\rho _{i+j}\rho _{i-j}$ for all $i,j\in Z.$

The representations above induce vector bundles over the space $BD_{4k}.$We
will denote the induced bundles over $BD_{4k}$ again by the same letters and
set the reductions $v_{1}=\eta _{1}-1,$ $v_{2}=\eta _{2}-1,$ $v_{3}=\eta
_{3}-1$ and $\phi =\rho _{1}-2$ in $K(BD_{4k})$.

As we said before, the results of [4] can be copied here with some
modifications. The relations above are exactly the same relations in [4]
written for the generalized quaternion group $Q_{4k}$, except that the start
relation and the end relation are interchanged! Because of that, it is
meaningful to choose the pair of representations $\eta _{2},\eta _{3}$
instead of the pair $\eta _{1},\eta _{2}$ in our problem as the generators
in the description. Note that, since $\eta _{3}=\eta _{1}\eta _{2},$ we have 
$v_{3}=v_{1}+v_{2}+v_{1}v_{2}$ and although highly convoluted, we can write
the description by using $v_{1}$ and $v_{2},$ whenever we want.

Now, since $\eta _{2},\eta _{3}$ and $\rho _{1}$ generate the representation
ring, $v_{2},v_{3}$ and $\phi $ generate $K(BD_{4k})$ due to ASCT. Now, we
can transform the relations in the representation ring to our new variables
and we can write the following:

\textbf{Theorem 2: }Let $n=2k$ and $k\geq 2.$ Then $K(BD_{2n})$ is generated
by $v_{2},v_{3}$ and $\phi $ and the minimal set of relations

\begin{center}
\begin{tabular}{ll}
$v_{2}^{2}=-2v_{2}$ & \textbf{(}$\text{\textbf{1)} }$ \\ 
$v_{3}^{2}=-2v_{3}$ & ($\text{\textbf{2)} }$ \\ 
$v_{2}\phi =\psi ^{k-1}(\phi )-\phi -2v_{2}-v_{3}\text{ }$ & \textbf{(3, }$%
\text{\textbf{for }}k\geq 3$ \textbf{and} $k$ \textbf{is odd)} \\ 
$v_{2}\phi =\psi ^{k-1}(\phi )-\phi -2v_{2}$ & \textbf{(3, }$\text{\textbf{%
for }}k\geq 3$ \textbf{and} $k$ \textbf{is even)} \\ 
$v_{3}\phi =-2v_{3}$ & \textbf{(4)} \\ 
$v_{2}v_{3}=4\phi +\phi ^{2}-2v_{2}-2v_{3}\text{ }$ & \textbf{(}$\text{%
\textbf{5, for }}k=2$\textbf{)} \\ 
$v_{2}v_{3}=\psi ^{k}(\phi )-2v_{2}-v_{3}$ & $\text{(\textbf{5, for }}k\geq
3 $ \textbf{and} $k$ \textbf{is odd)} \\ 
$v_{2}v_{3}=\psi ^{k}(\phi )-2v_{2}\text{ }$ & \textbf{(}$\text{\textbf{5,
for }}k\geq 3$ \textbf{and} $k$ \textbf{is even)}%
\end{tabular}
\end{center}

\textit{Proof: }Thanks to [4], we derived these relations quickly from the
Theorem 1, [4]. But, there is a slight modification for the case $k$ is odd
in the relations 3 and 5 above. Now, let us analyze this quick and
surprising result through the AHSS, confirming the minimality of the
relations and thus completing the proof of the theorem. Before that, let us
recall the cohomology.

The integral cohomology of the group $D_{2n}$ is given by, [2],

\begin{equation*}
H^{p}(BD_{2n};Z)=\left\{ 
\begin{tabular}{l}
$Z$ \\ 
$Z_{n}\oplus (Z_{2})^{2s}$ \\ 
$(Z_{2})^{2s}$ \\ 
$(Z_{2})^{2s+2}$ \\ 
$(Z_{2})^{2s+1}$%
\end{tabular}%
\right. 
\begin{tabular}{l}
$p=0$ \\ 
if $p=4s,$ $s\geq 1$ \\ 
if $p=4s+1$ \\ 
if $p=4s+2$ \\ 
if $p=4s+3$%
\end{tabular}%
\end{equation*}%
and we can read the filtrations of the second page of the AHSS from the
cohomology described above. But, since, the odd dimensional cohomology is
not trivial, AHSS may not collapse on page 2. And in fact, this is the case
here.

It is enough to consider $k=2$ in order to explain the minimilatiy of the
relations through the AHSS, since exactly the same explanation goes for $%
k\geq 3$.

Now, we have two natural group inclusions $Z_{2}<D_{8}$ and they induce two
natural ring homomorphisms $j_{1}:K(BD_{8})\rightarrow K(BZ_{2})$ and $%
j_{2}:K(BD_{8})\rightarrow K(BZ_{2}).$ Recall that $%
K(BZ_{2})=Z[v]/(v^{2}+2v).$ Now, since under $j_{1}$ and $j_{2},$ $v_{1}$
and $v_{2}$ maps to $v,$ we deduce that $E_{\infty }^{2,-2}=Z_{2}\oplus
Z_{2} $ of the AHSS of $K(BD_{8})$ is generated by $v_{1}$ and $v_{2}.$ Note
that $v_{1}=v_{2}+v_{3}+v_{2}v_{3}$ here.

For $K(BD_{8}),$ the main relation is $4\phi
=-v_{2}^{2}-v_{3}^{2}+v_{2}v_{3}-\phi ^{2}.$ From this, we have $4\phi =$ $%
h.o.t.$ and we claim that the $Z_{4}$ part of the filtration $E_{\infty
}^{4,-4}=Z_{4}\oplus Z_{2}\oplus Z_{2}$ is generated by $\phi $. Really, we
have a natural group inclusion $Z_{4}<D_{8},$ and this induces a natural
ring homomorphisms $j:K(BD_{8})\rightarrow K(BZ_{4})$. Recall that $%
K(BZ_{4})=Z[\mu ]/((1+\mu )^{4}-1)$. Under $j,$ $\phi $ maps to the infinite
alternating sum $\mu ^{2}-\mu ^{3}+\mu ^{4}-...$ and since $\mu ^{2}$
generates $E_{\infty }^{4,-4}$ filtration of $K(BZ_{4}),$ by naturailty, $%
\phi $ should generate a direct summand of $E_{\infty }^{4,-4}$ for $%
K(BD_{8})$ too. And there is one choice.

We claim that $Z_{2}\oplus Z_{2}$ part of $E_{\infty }^{4,-4}$ is also
surviving and is generated by $v_{1}v_{3}$ and $v_{2}v_{3}$ respectively.
First of all, $j_{1}$ and $j_{2}$ should be non-trivial on these
filtrations, because $E_{\infty }^{4,-4}$of $K(BZ_{2})$'s are generated by $%
v^{2}$. Since, under $j_{1}$ and $j_{2},$ $\phi $ which generates $Z_{4}$
part of the filtration $E_{\infty }^{4,-4},$ maps to zero, only the other
parts can provide this non-triviality. On the other hand, under $j_{1}$ and $%
j_{2}$, the element $v_{3}$ goes to $v$ so that $v_{1}v_{3}$ and $v_{2}v_{3}$
both go to $v^{2}$. And the claim is proved. Note that, $Z_{2}\oplus Z_{2}$
part of $E_{\infty }^{4,-4}$ can not be generated by $v_{1}^{2}$ and/or $%
v_{2}^{2}.$ This is very confusing, but, under $j,$ since both $v_{1}$ and $%
v_{2}$ go to $2\mu +\mu ^{2},$ $v_{1}^{2}$ and $v_{2}^{2}$ go to $-2\mu
^{3}-3\mu ^{4}-\mu ^{5},$ and they are related to the filtration $E_{\infty
}^{6,-6}.$

So far, we noticed that the AHSS is collapsing on the second page for the
filtrations $E_{\infty }^{2,-2}$ and $E_{\infty }^{4,-4}.$ But, this will
not be the case for the higher filtrations.

We claim that $E_{\infty }^{6,-6}$ is $Z_{2}\oplus Z_{2}\oplus Z_{2}$ and
that it is generated by $v_{1}v_{3}^{2}$, $v_{2}v_{3}^{2}$ and $v_{1}v_{2}.$
Note that on the second page, $E_{2}^{6,-6}$ is $Z_{2}\oplus Z_{2}\oplus
Z_{2}\oplus Z_{2}$ and it has four summands. We claim that one summand $%
Z_{2} $ can not survive to infinity. This claim is not difficult to prove
after the observations we did on the filtrations $E_{\infty }^{2,-2}$ and $%
E_{\infty }^{4,-4}.$ First of all, $v_{1}v_{3}^{2}$, $v_{2}v_{3}^{2}$ should
survive on the filtration $E_{\infty }^{6,-6},$ because they go to $v^{2}$
under $j_{1}$ and $j_{2}.$ This takes care of two summands. The product $%
v_{1}v_{2}$ has clearly on the filtration $E_{\infty }^{6,-6}$ because under 
$j$ it goes to where $v_{2}^{2}$ goes, due to $%
v_{2}^{2}=v_{1}v_{2}+v_{2}v_{3}.$ And this takes care of the third summand.
Because of the relations, we exhausted all elements that can live on the
filtration $E_{\infty }^{6,-6}$ and we are done.

Now, the picture is clear for $K(BD_{8})$ for the higher filtrations:
Firstly, the filtration $E_{\infty }^{4t,-4t}$ is $Z_{4}\oplus Z_{2}\oplus
Z_{2}$ and it is generated by $\phi $ and $v_{1}v_{3}^{2t-1}$, $%
v_{2}v_{3}^{2t-1}$ respectively, where $t\geq 1$. Secondly, the filtration $%
E_{\infty }^{4t+2,-4t-2}$ is $Z_{2}\oplus Z_{2}\oplus Z_{2}$ and it is
generated by $v_{1}v_{3}^{2t}$, $v_{2}v_{3}^{2t}$ and $v_{1}v_{2}^{t}$
respectively, where $t\geq 1$.

Finally, $2v_{2}$ can be expressed by $\phi $, $v_{2}v_{3}$ for the even
case and by $\phi $, $v_{2}v_{3}$ and $v_{3}$ for the odd case by means of
the relation 5. Since $v_{3}=v_{1}+v_{2}+v_{1}v_{2},$ we deduce that $2v_{2}$
can be expressible by the elements we chose as the generators of the
filtrations of the AHSS.$\blacksquare $

We also note that the generator $\phi $ satisfies the relation $g_{2k}(\phi
)=\psi ^{k+1}(\phi )-\psi ^{k-1}(\phi )=0$. The polynomial $g_{2k}(\phi )$
is explicitly given in [4]. It should be clear to the reader that this
polynomial is the even version of the magical polynomial $wf_{n}(w)$ of the
odd case.

As a final note, we point out the very natural continuation of the problem: $%
K$-rings of the classifying spaces of the symmetric groups. For $n=3,$ $%
D_{6}=S_{3}$ and $K(BS_{3})$ is explained here. But, for $n\geq 4,$ the ring 
$K(BS_{n})$ becomes devilishly complicated. The polynomials arising as
minimal relations of these rings can have tremendous connections with the
homotopy groups of spheres.

\end{document}